\documentclass{amsart}
\title{An Approach to Hopf Algebras via Frobenius Coordinates I}
\author{Lars Kadison and A.A. Stolin}
\address{Chalmers University of Technology/G{\" o}teborg
University \\
S-412 96 G{\" o}teborg, Sweden}
\email{kadison@math.ntnu.no \\
astolin@math.chalmers.se}
\date{}
\newtheorem{theorem}{Theorem}[section]
\newtheorem{lemma}[theorem]{Lemma}
\newtheorem{prop}[theorem]{Proposition}
\newtheorem{corollary}[theorem]{Corollary}

\theoremstyle{definition}
\newtheorem{definition}[theorem]{Definition}

\newtheorem{remark}[theorem]{Remark}

\newcommand{\End}{\fam0\rm End}
\newcommand{\Hom}{\fam0\rm Hom}

\newcommand{\Id}{\fam0\rm Id}

\newcommand{\forevery}{\fam0\forall \,}
\begin{document}
\begin{abstract}
In Section 1 we  introduce Frobenius coordinates in the general setting that includes  
 Hopf subalgebras.   
In Sections 2 and 3 we review briefly the theories 
of Frobenius algebras   and augmented Frobenius algebras with some new material
in Section 3. In Section 4 we study the Frobenius
structure of an FH-algebra $H$ \cite{Par72}
and extend two recent theorems in \cite{EG}.  We obtain two Radford formulas 
for the antipode in $H$ and generalize in Section 7 the results on its order in \cite{FMS}.
We study the Frobenius structure on an FH-subalgebra pair in Sections 5 and 6.  In Section 8 
we show that the quantum double of $H$ is symmetric and unimodular.  
\end{abstract}   
\maketitle

\section{Introduction}

 Suppose $A$ and $S$ are noncommutative associative rings with $S$ a 
unital subring in $A$, or 
stated equivalently, $A/S$ is a ring extension. Given a  ring automorphism
$\beta:S \rightarrow S$, a left $S$-module $M$
receives the $\beta$-twisted module structure ${}_{\beta}M$
by $s \cdot_{\beta} m := \beta(s)m$ for each $s \in S$ and $m \in M$.
  $A/S$ 
 is said to be a $\beta$-Frobenius extension if 
the natural module $A_S$ is finite projective,
  and 
$$ {}_SA_A \cong {}_{\beta}\Hom_S(A_S,S_S)_A$$  \cite{FMS,Par64}.  A very 
useful characterization of $\beta$-Frobenius extensions is that 
they are the ring extensions
having a Frobenius coordinate system.   
 A {\it Frobenius coordinate system} for 
 a ring extension $A/S$ is 
 data $( E,x_i,y_i )$ where 
 $E: {}_SA_S \rightarrow {}_{\beta}S_S$ is a bimodule homomorphism, 
called the {\it Frobenius homomorphism}, and elements 
$x_i,y_i \in A$ $(i = 1,\ldots,n)$, called  {\it dual bases},  such that
 for every $a \in A$:
\begin{equation}
\sum_{i=1}^n \beta^{-1}(E(ax_i)) y_i = a = \sum_i x_iE(y_i a).
\label{eq:dualbase}
\end{equation}

One of the most important points  about 
Frobenius coordinates for $A/S$ is that any two of these, 
$(E,x_i,y_i)$ and $(F,z_j,w_j)$,    
differ by only an invertible $d \in C_A(S)$, the
centralizer of $S$ in $A$:  viz.\ $F = Ed$
and $\sum_i x_i \otimes d^{-1}y_i = \sum_j z_j \otimes w_j$ \cite{Par64}.    
 The {\it Nakayama automorphism} $\eta$
of $C_A(S)$ may be  defined by \[
E(\eta(c)a) = E(ac)
\] for every $a \in A, c \in C_A(S)$.  Then from  Equations~\ref{eq:dualbase},
$
\eta(c) = \sum_i \beta^{-1}(E(x_ic))y_i ,
$
and 
\begin{equation}
\eta^{-1}(c) = \sum_i x_i E(cy_i).
\label{eq:formue}
\end{equation} 
The  Nakayama automorphisms $\eta$ and $\gamma$ relative to two Frobenius homomorphisms $E$
and $F= Ed$, respectively,  
are related by $\eta \gamma^{-1}(x) = dxd^{-1}$ for every $x \in C_A(S)$ 
\cite{Par64}. If $A$ is a $k$-algebra and $S = k1_A$, then $\beta$
is necessarily the identity by a short calculation \cite{NT} and 
 $C_A(S) = A$.  

For example, a Hopf subalgebra $K$ in a finite dimensional Hopf algebra $H$
over a field is a free $\beta$-Frobenius extension.  The natural
module $H_K$ is free by the  theorem of Nichols-Zoeller \cite{NZ}. 
  By a theorem of Larson-Sweedler in \cite{LS}, the antipode
is bijective, and $H$ and $K$
are Frobenius algebras with Frobenius homomorphisms which are left
or right integrals in the dual algebra. From
Oberst-Schneider \cite[Satz 3.2]{OS73} we have a formula (cf.\ Equation~\ref{eq:soon})
 that implies that the Nakayama automorphism of $H$, $\eta_H$, 
restricts to a mapping of $K \rightarrow K$. 
It follows from
Pareigis \cite[Satz 6]{Par64} that $H/K$ is a $\beta$-Frobenius extension, where 
the automorphism
$\beta$ of $K$ is the following
 composition of the Nakayama automorphisms of $H$ and
$K$:
\begin{equation}
\beta = \eta_K \circ \eta_H^{-1}
\label{eq:beta}
\end{equation}
(cf.\ Section~5). 

This paper continues our investigations
in \cite{BFS,BFS2,LK96,NEFE} on the interactions of Frobenius
algebras/extensions with Hopf algebras. We apply
  Frobenius coordinates to a class
of Hopf algebras over commutative rings called FH-algebras,
which are Hopf algebras that are simultaneously Frobenius algebras
(cf.\ Section~4). This
class was introduced in \cite{Par71, Par72} and   
includes the 
finite dimensional Hopf algebras as well as the 
 finite projective Hopf algebras over commutative rings with 
trivial Picard group
(such as semi-local or polynomial rings).
The added generality would apply for example to a  
Hopf algebra $H$ over  a Dedekind
domain $k$ satisfying the condition that the element represented by 
the $k$-module of left integrals $\int^{\ell}_{H^*}$ in the Picard
group of $k$ be trivial.

This paper is organized as follows. 
In Section 2, we  review the 
basics of Frobenius algebras and Frobenius coordinates, as well
as separability. In Section 3, we study  
 norms, integrals and modular functions for augmented Frobenius
algebras over a commutative ring, giving a lemma on the effect of automorphisms and anti-automorphism
on s. In Section 4,
we derive by means of different Frobenius coordinates Radford's
Formula~\ref{eq:friday} 
for $S^4$ and Formula~\ref{eq:Radford} relating $S^2$, $t_1, t_2$, 
where $t$ is a right norm for $H$. 
This extends two formulas
 in \cite[$k$ = field]{Rad76,Rad94} to  FH-algebras
with  different proofs. Then we generalize two recent results
 of Etingof and Gelaki \cite{EG}, the main one stating
 that a finite dimensional semisimple and cosemisimple
 Hopf algebra   is involutive. We show that
with a small condition on $2 \in k$ a separable and coseparable Hopf $k$-algebra  is 
 involutive (Theorem~\ref{th-geneg}).  Furthermore,
if $H$ is separable and satisfies a certain bound on its local ranks,
then $H$ is coseparable and therefore involutive (Theorem~\ref{th-saturn}).  
 
In Section 5, we prove that a subalgebra pair of FH-algebras $H \supseteq K$
is a  $\beta$-Frobenius extension, though not necessarily free.
In Section 6, we derive by means of different Frobenius 
coordinates  Equation~\ref{eq:cecilia}
relating the different elements in a $\beta$-Frobenius coordinate system for a Hopf
subalgebra pair $K \subset H$ given by
Fischman-Montgomery-Schneider \cite{FMS} . 
In Section 7, we prove 
 that a group-like element in a finite projective Hopf algebra $H$ 
over a Noetherian ring $k$
has finite order dividing the least common multiple $N$
of the $P$-ranks of $H$ as a $k$-module. From the theorems 
in Section 4 it follows 
that $S$ has  order dividing $4N$, and, should $H$ be an FH-algebra,
that the Nakayama automorphism $\eta$
has finite order dividing $2N$, as obtained 
for fields in \cite{Rad76} and \cite{FMS}, respectively. 
In Section 8, we  extend 
the Drinfel'd notion of quantum double to FH-algebras, then prove
that the quantum double of an 
FH-algebra $H$  is a unimodular and symmetric FH-algebra.

\section{A brief review of Frobenius algebras}

All rings in this
paper have
$1$, homomorphisms preserve $1$, and unless otherwise specified
 $k$ denotes a commutative ring. 
Given an associative, unital $k$-algebra $A$,  $A^*$ denotes the dual module $\Hom_k(A,k)$,
which is an $A$-$A$ bimodule as follows: given $f \in A^*$ and $a\in A$,
$af$ is defined by $(af)(b) = f(ba)$ for every $b \in A$, while $fa$ is defined by
$(fa)(b) = f(ab)$.  We also consider the {\it tensor-square},
$A \otimes A$ as a natural $A$-bimodule  given
by $a(b \otimes c) = ab \otimes c$ and $(a \otimes b)c :=
a \otimes bc$ for every $a,b,c \in A$. An element $\sum_i z_i \otimes w_i$
 in the tensor-square
is called {\it symmetric} if it is left fixed by the transpose map given
by $a \otimes b \mapsto b \otimes a$ for every $a,b \in A$. 

We first consider some preliminaries on a 
Frobenius algebra $A$ over a commutative ring $k$. 
$A$ is a {\it Frobenius algebra} if the natural module  
 $A_k$ is  {\it finite projective} (= finitely generated projective), and 
\begin{equation}
A_A \cong A^*_A.
\end{equation}
 Suppose $f_i \in A^*$,
$x_i \in A$ form a finite projective base, or dual bases, of $A$ over $k$:
i.e., for every $a \in A$, $\sum_i x_i f_i(a) = a$.
Then there are $y_i \in A$ and a cyclic generator
 $\phi \in A^*$ such that the $A$-module
isomorphism is given by $a \mapsto \phi a$, and
\begin{equation}
\sum_i x_i \phi(y_ia) = a = \sum \phi(ax_i)y_i,
\label{eq:june3}
\end{equation}
for all $a \in A$.  It follows that $\phi$ is  nondegenerate
(or faithful) in the
following sense:  a linear functional $\phi$ on an algebra $A$ 
is {\it nondegenerate} if $a,b \in A$ such that 
$a \phi = b \phi$ or  $\phi a = \phi b$ 
implies $a = b$.

We refer to $\phi$ as a {\it Frobenius homomorphism},
$(x_i,y_i)$ as  {\it dual bases}, and $(\phi,x_i,y_i)$ as a {\it 
Frobenius system} or \textit{Frobenius coordinates}.  It is useful to note from the start that $xy = 1$
implies $yx = 1$ in $A$, since an epimorphism of $A$ onto itself is
automatically bijective \cite{Par71,Sil}. 

It is equivalent to define a $k$-algebra $A$ Frobenius if 
 $A_k$ is finite projective and ${}_AA \cong {}_AA^*$.  In fact, with $\phi$ 
defined above, the mapping $a \mapsto a\phi$ is such an isomorphism,
by an application of Equations~\ref{eq:june3}. 

Note that the bilinear form on $A$ defined by $\langle a,b \rangle
:= \phi(ab)$ is a nondegenerate inner product 
which is associative: $\langle ab,c \rangle = \langle a,bc \rangle$ for
every $a,b,c \in A$. 

The Frobenius homomorphism  is unique up to an invertible element in $A$.  
If $\phi$ and $\psi$ are  Frobenius homomorphisms for $A$, then $\psi = d\phi$
for some $d \in A$.  Similarly, $\phi = d' \psi$ for some $d' \in A$,
from which it follows that $dd' = 1$.  The element $d$ is referred to as 
the (left) {\it derivative} $\frac{d\psi}{d\phi}$ of $\psi$ with respect to $\phi$.
  Right derivatives in the group of units $A^{\circ}$ of $A$
 are similarly defined.

If $(\phi,x_i,y_i)$ is a Frobenius system for $A$, then 
$e := \sum_i x_i \otimes y_i$
is an  element in the tensor-square $A \otimes_k A$
which is independent of the choice of dual bases for $\phi$,
 called the {\it Frobenius element}.
 By a computation 
involving Equations~\ref{eq:june3},
 $e$ is a 
Casimir element
satisfying $ae = ea$ for every $a \in A$, whence $\sum_i x_iy_i$ is in the center
of $A$.\footnote{$e$ is the 
transpose of the element $Q$ in \cite{BFS2}.}  
 It follows that $A$ is $k$-separable if and only if
there is a $a \in A$ such that 
\begin{equation} 
\sum_i x_i a y_i = 1.
\label{eq:bursdag}
\end{equation}

For each $d \in A^{\circ}$, we easily check that $(\phi d, x_i, d^{-1}y_i)$ and
$(d\phi,x_i d^{-1},y_i)$ are the other Frobenius systems in 
a one-to-one correspondence. It follows that a Frobenius element is also
unique, up to a unit in $A \otimes A$ (either $1 \otimes d^{\pm 1}$ or
$d^{\pm 1} \otimes 1$).

A {\it symmetric algebra} is a Frobenius algebra $A/k$
 which satisfies the stronger condition: 
\begin{equation}
{}_AA_A \cong {}_A(A^*)_A.
\label{eq:symmetry}
\end{equation}
Choosing an isomorphism $\Phi$, the linear functional $\phi := 
\Phi(1)$ is a Frobenius homomorphism satisfying
\(
\phi(ab)  = \phi(ba)
\)
for every $a, b \in A$: i.e., $\phi$ is an  trace on $A$. 
The dual bases $x_i,y_i$ for this $\phi$ forms a symmetric element
in the tensor-square, since for every $a \in A$, 
\begin{eqnarray}
\sum_i ax_i \otimes y_i & = & \sum_{i,j} y_j  \otimes \phi(ax_i x_j)y_i \nonumber \\
& = & \sum_j y_j \otimes x_j a . 
\label{eq:easter}
\end{eqnarray}

A $k$-algebra $A$ with $\phi \in A^*$ and $x_i,y_i \in A$
satisfying either $\sum_i x_i \phi(y_i a) = a$ for every $a \in A$
or $\sum_i \phi(ax_i)y_i = a$ for every $a \in A$ is automatically Frobenius.
As a corollary, one of the dual bases equations implies the other. 
For if $\sum_{i=1}^n (x_i\phi) y_i = \Id_A$, then $A$ is
explicitly finite projective over $k$, and it follows that $A^*$ is
finite projective too.  The homomorphism ${}_AA \rightarrow
{}_AA^*$ defined by $a \mapsto a\phi$ for all $a \in A$
 is surjective, since
given $f \in A^*$, we note that $f = (\sum_i f(y_i) x_i)\phi$.
Since $A$ and $A^*$ have the same $P$-rank for each
prime ideal $P$ in $k$, the epimorphism $a \mapsto a\phi$ is bijective
\cite{Sil}, whence  ${}_AA \cong {}_AA^*$.
Starting with the other equation in the hypothesis, we similarly
prove that $a \mapsto \phi a$ is an isomorphism $A_A \cong A^*_A$. 

The Nakayama automorphism of a Frobenius algebra $A$ is an algebra 
automorphism $\alpha: A \rightarrow A$ defined
by 
\begin{equation}
  \phi \alpha(a) = a\phi
\label{eq:1}
\end{equation}
  for every $ a \in A$. In terms of the associative inner product,
$\langle x,a \rangle = \langle \alpha(a),x \rangle$ for every $a,x \in A$.
 $\alpha$  is an inner automorphism iff
$A$ is a symmetric algebra. The Nakayama automorphism $\eta$ of another Frobenius
homomorphism $\psi = \phi d$, where $d \in A^{\circ}$, is
given by \begin{equation}
\eta(x) = \sum_i \phi(dx_ix)d^{-1}y_i = \sum_i d^{-1} \phi(\alpha(x)dx_i)y_i
= d^{-1}\alpha(x) d,
\label{eq:alex} 
\end{equation} 
so that $\alpha \eta^{-1}(x) = dxd^{-1}$.  Thus
 the Nakayama automorphism is unique up to an inner automorphism. 
A Frobenius algebra $A$ is a symmetric algebra if and only if its
Nakayama automorphism is inner. 

Another formula for $\alpha$ is obtained from Equations~\ref{eq:1}
and~\ref{eq:june3}: for every $a \in A$, 
\begin{equation}
\alpha(a) = \sum_i \phi(x_i a) y_i.
\end{equation}
If the Frobenius element $\sum_i x_i \otimes y_i$ is symmetric, it
follows from this equation that $\alpha = \Id_A$.  Together with
Equation~\ref{eq:easter}, this proves:

\begin{prop}
A Frobenius algebra $A$ is a symmetric algebra if and only if
it has a symmetric Frobenius element. 
\end{prop}

Equation~\ref{eq:symmetry} generalizes
to all Frobenius algebras as follows.
A Frobenius isomorphism $\Psi: A_A \stackrel{\cong}{\rightarrow} A^*_A$ 
induces a bimodule isomorphism
where one bimodule is twisted by the Nakayama automorphism $\alpha$:
\begin{equation}
{}_AA_A \cong {}_{\alpha^{-1}}A^*_A, 
\label{eq:yama}
\end{equation}
since with $\phi = \Psi(1)$ Equation~\ref{eq:1} yields 
\[
\Psi(a_1aa_2) = \phi a_1 aa_2
= \alpha^{-1}(a_1)\phi a a_2 = \alpha^{-1}(a_1) \Psi(a)a_2. \]

The left and right derivatives of a pair of Frobenius homomorphisms
differ by an application of the Nakayama automorphism (cf.\ Equation~\ref{eq:1}).
A computation applying  Equations~\ref{eq:june3} and~\ref{eq:1} 
proves that for every $a \in A$,
\begin{equation} 
\sum_i x_i a \otimes y_i = \sum_i x_i \otimes \alpha(a) y_i .
\label{eq:laurdag}
\end{equation}

In closing this section, we refer the reader to \cite{EN,BFS,NEFE}
for more on Frobenius algebras over commutative
rings, and to \cite{Yam} for a survey of the representation theory of Frobenius over fields
and work on the Nakayama conjecture. 

\section{Augmented Frobenius algebras}

A $k$-algebra  $A$ is said to be an {\it augmented algebra}
if there is an algebra homomorphism $\epsilon: A \rightarrow k$,
called an {\it augmentation}. 
An element  $t \in A$ satisfying
$ta = \epsilon(a)t
$,  $\forevery a \in A$,
is called a {\em right integral} of $A$. 
  It is clear that the set of right
integrals, denoted by $\int^r_A$, is a two-sided ideal of $A$, since
for each $a \in A$, the element $at$ is also a right
integral.  Similarly for the space of left integrals, denoted by
$\int_A^{\ell}$. If $\int^r_A = \int^{\ell}_A$, 
$A$ is said to be {\it unimodular}. 

Now suppose that $A$ is a Frobenius algebra with augmentation $\epsilon$. 
We claim that a nontrivial right integral exists in  $A$.
Since $A^* \cong A$ as right $A$-modules, 
an element   $n \in A$ exists such that  
$\phi n  = \epsilon$ where $\phi$ is a Frobenius homomorphism.
Call $n$ the {\it right norm} in $A$ with respect to $\phi$. 
 Given $a \in A$, we compute in $A^*$: 
 \[
\phi na =(\phi n)a = \epsilon a =\epsilon(a)\epsilon =
 \phi n\epsilon(a). \]
By nondegeneracy of $\phi$, 
 $n$ satisfies $na = n \epsilon(a)$ for every $a \in A$.

\begin{prop}  
If $A$ is an augmented Frobenius algebra, then the set $\int^r_A$ of right
integrals is a two-sided ideal which is free cyclic $k$-summand of $A$ 
generated by a right norm. 
\label{prop-ein}
\end{prop}
\begin{proof} 
The proof is based on \cite[Theorem 3]{Par71}, which assumes that $A$
 is also a Hopf algebra. 
 Let  $\phi \in A^*$ be a  Frobenius homomorphism, and 
 $n \in A$ satisfy $\phi n = \epsilon$, the augmentation.  
Given a  right integral $t \neq 0$, we note that 
\[ \phi t = \phi(t) \epsilon = \phi(t) \phi n = \phi n\phi(t),\]
whence 
\begin{equation}
t = \phi(t)n.
\label{eq:epif}
\end{equation}
  Then $\langle n \rangle := \{ \rho n |\, \rho \in k \}$
coincides with the set of all right integrals. 

Given $\lambda \in k$ such that $\lambda n = 0$, it follows that
\[
\phi(n) \lambda = \epsilon(1)\lambda = \lambda = 0, \]
whence $\langle n \rangle$ is a free $k$-module.  
Moreover, $\langle n \rangle$ is a direct $k$-summand in $A$ since
$a \mapsto \phi(a)n$ defines a $k$-linear projection of $A$ onto $\langle n \rangle$.  
 \end{proof}

The right norm in $A$ is unique up to a unit in $k$, since
norms are free generators of $\int^r_A$ by the proposition. 
The notions of norm and integral only coincide if $k$ is a field. 

Similarly the space $\int^{\ell}_A$ of left integrals is a rank one free summand
 in $A$, generated by any left norm.  In general 
 the spaces of right and left integrals do not coincide, and one defines
an augmentation on $A$ that measures the deviation from unimodularity. 
 In the notation of the
proposition and its proof, for every $a \in A$, the element $an$ is
a right integral since the right norm $n$ is.  From
Equation~\ref{eq:epif} one concludes that $an = \phi(an)n = (n\phi)(a)n$.
The function 
\begin{equation}
m := n\phi: A \rightarrow k
\end{equation}
 is called the {\it right modular function}, which is 
an augmentation since $\forevery a,b \in A$ we have
$(ab)n = m(ab)n = a(bn) = m(a)m(b)n$ and $ n $  is a  free generator of $\int^r_A$.

The next proposition and corollary we believe
has not been noted in the literature before. 

\begin{prop}
If $A$ is an augmented Frobenius algebra 
and $\alpha$ the Nakayama automorphism,
then in the notation  above, 
\begin{equation}
m \circ \alpha = \epsilon.
\label{eq:onsdag}
\end{equation}
\end{prop}
\begin{proof}
We note that $\phi \circ \alpha = \phi$ by evaluating each side of Equation~\ref{eq:1}
on $1$. Then for each $x \in A$, 
\[
m(\alpha(x)) = (n\phi)(\alpha(x)) = (\phi \alpha(n))(\alpha(x)) = (\phi \circ \alpha) (x) = \phi(x). \qed
\]
\renewcommand{\qed}{}\end{proof}

The next corollary follows from noting that if $\alpha$ is an inner automorphism,
then $m = \epsilon$ from the proposition. 
\begin{corollary}
If $A$ is an augmented symmetric algebra, then $A$ is unimodular.
\label{cor-interesting}  
\end{corollary} 

We note two useful identities for the right norm,  
\begin{eqnarray}
n & =&  \sum_i \phi(nx_i)y_i = \sum_i \epsilon(x_i)y_i 
\label{eq:maritshoes}\\
& = & \sum_i x_i (n\phi)(y_i) = \sum_i x_i m(y_i).
\label{eq:nesbyen}
\end{eqnarray}

As an example, consider  $A := k[X]/(X^n)$  where $k$ is a commutative ring
and $aX = X\alpha(a)$ for some automorphism $\alpha$ of $k$
and every $a \in k$. Then $A$ is an augmented Frobenius algebra
with Frobenius homomorphism $\phi(a_0 + a_1X + \cdots + a_{n-1}X^{n-1}) :=
a_{n-1}$, dual bases $x_i = X^{i-1}$, $y_i = X^{n-i}$ ($i = 1,\cdots,n$),
and augmentation $\epsilon( a_0 + a_1X + \cdots + a_{n-1}X^{n-1}) := 
a_0$.  It follows that a left and right norm is given by $n = X^{n-1}$, and 
$A$ is   symmetric and unimodular.  $A$ is not a Hopf algebra unless $n$ is
a prime $p$ and the characteristic of $k$ is $p$ (cf.\ \cite{FMS}).

 The next proposition is well-known for finite dimensional Hopf algebras
\cite{Mont}. 
\begin{prop}
Suppose $A$ is a separable augmented Frobenius algebra.
Then $A$ is unimodular.
\label{prop-marchsnow}
\end{prop}
\begin{proof}
The Endo-Watanabe theorem in  \cite{EW} 
states that separable projective algebras are symmetric algebras.
The result  follows then from Corollary~\ref{cor-interesting}. 
\end{proof}

 We will use   repeatedly in Section 4 
several general principles summarized in the next lemma.
Items 1, 2 and 3 below are valid without the assumption of augmentation or $\epsilon$-invariance. 

\begin{lemma}
  Suppose  $(A,\epsilon)$ is an augmented Frobenius algebra
and  $\alpha$ (respectively, $\beta$) is a $k$-algebra automorphism (resp.\ anti-automorphism)
of $A$ satisfying  $\epsilon$-invariance: viz.\ 
 $\epsilon \circ \alpha = \epsilon$. Let  $(\phi_A,x_i,y_i)$ be 
Frobenius coordinates of $A$.  Then
\begin{enumerate}
\item The Frobenius system is transformed by 
$\alpha$ into a Frobenius system
\[
(\phi_A \circ \alpha^{-1}, \alpha(x_i), \alpha(y_i)).\]
\item The Frobenius system  is transformed by $\beta$
 into the Frobenius system
\[
(\phi_A \circ \beta^{-1}, \beta(y_i), \beta(x_i)). \]
\item If $B$ is another Frobenius $k$-algebra with Frobenius homomorphism
 $\phi_B$, then $A \otimes B$ is a Frobenius algebra with
Frobenius homomorphism $\phi_A \otimes \phi_B: A \otimes B \rightarrow k$.
\item  $\alpha$ sends integrals to integrals and norms to norms, respecting
chirality. 
\item $\beta$ sends integrals to integrals and norms to norms, reversing
chirality.
\end{enumerate}
\label{lemma-emma}
\end{lemma}
\begin{proof}  1 is proven by applying $\alpha$ to $\sum_i \phi_A(ax_i)y_i = a$,
obtaining \[
\sum_i \phi_A \alpha^{-1}(\alpha(a)\alpha(x_i))\alpha(y_i) = \alpha(a)
\]
for every $a \in A$.  2 is proven similarly. 3 is easy. 
4 is proven by applying
first $\alpha$ to $ta = \epsilon(a)t$, obtaining that $\alpha(t) \in \int^r_A$
if $t$ is too. Next, if $\phi_A n = \epsilon$, then $(\phi_A \circ \alpha^{-1})\alpha(n) =
\epsilon$ as well, which together with 1 proves 4. 5 is proven similarly. 
\end{proof} 

\section{FH-algebras}

We continue with $k$ as a commutative ring. 
We review the basics of a Hopf algebra $H$ which is finite projective over
 $k$ \cite{Par71}. A bialgebra $H$ is an algebra and coalgebra
where the comultiplication and the counit are algebra homomorphisms. 
 We use a reduced Sweedler notation given by 
  \[
\Delta(a) = \sum_{(a)} a_{(1)} \otimes a_{(2)} := \sum a_1 \otimes a_2
\] 
for the values of the comultiplication homomorphism $H \rightarrow 
H \otimes_k H$.
The counit is the $k$-algebra homomorphism
$\epsilon: H \rightarrow k$ and satisfies $\sum_i \epsilon(a_1)a_2 =
\sum a_1 \epsilon(a_2) = a$ for every $a \in H$.

A Hopf algebra $H$ is a bialgebra with antipode.  
The antipode $S: H \rightarrow H$ is an anti-homomorphism of algebras and
coalgebras satisfying $\sum S(a_1)a_2 = \epsilon(a)1 = \sum a_1 S(a_2)$
for every $a \in H$.  

A {\it group-like element} in $H$ is defined
to be a $g \in H$ such that $\Delta(g) = g \otimes g$ {\it and} $\epsilon(g) =
1$. It follows that $g \in H^{\circ}$ and $S(g) = g^{-1}$. 

Finite projective Hopf algebras   enjoy the  duality
properties of finite dimensional Hopf algebras. 
 $H^*$ is  a Hopf algebra with convolution
product $(fg)(x) := \sum f(x_1)g(x_2)$.
The counit is given by $f \mapsto f(1)$. 
The unit of $H^*$ is the counit of $H$. The comultiplication on
$H^*$ is given by $\sum f_1 \otimes f_2 (a \otimes b) = f(ab)$ for every
$f \in H^*$, $a,b \in H$. The antipode is the dual of $S$,
a mapping of $H^*$ into $H^*$, denoted again by $S$
when the context is clear. 
Note that an augmentation $f$ in $H^*$ is a group-like element in $H^*$,
and {\it vice versa}, with inverse given
by $Sf = f \circ S$.

  As Hopf algebras,
$H \cong H^{**}$, the isomorphism being given by $x \mapsto {\rm ev}_x$,
the evaluation map at $x$:  
we fix this isomorphism as an identification of $H$ with $H^{**}$.  
The usual  left and right action of
an algebra  on its dual specialize to  the 
left action  of $H^*$
on $H^{**} \cong H$  given  by $g \rightharpoonup a := \sum a_1 g(a_2)$,
and the right action given by
$a \leftharpoonup g := \sum g(a_1)a_2$. 

We recall the definition of an equivalent version of Pareigis's  FH-algebras \cite{Par72}. 
\begin{definition}
A $k$-algebra $H$ is an FH-algebra if $H$ is a bialgebra and a Frobenius
algebra with Frobenius homomorphism $f$ a right integral in $H^*$.
Call $f$ the FH-homomorphism.\footnote{The authors have called
FH-algebras Hopf-Frobenius algebras in an earlier preprint.}
\label{def-june3}
\end{definition}

The condition that $f \in \int^r_{H^*}$ is equivalent to 
\begin{equation}
\sum f(a_1)a_2 = f(a)1
\end{equation}
for every $a \in H$.  
Note that $H$ is an augmented
Frobenius algebra with augmentation $\epsilon$. Let
 $t \in H$ be a right norm such that $ft = \epsilon$. Note that $f(t) = 1$. 
Fix the notation $f$ and $t$ for an FH-algebra. 
We show below that an FH-homomorphism is unique up to an invertible scalar in $k$. 
If $H$ is an FH-algebra and a symmetric algebra, we say that $H$ is a {\it
symmetric FH-algebra}. 

It follows from \cite[Theorem 2]{Par71} that an FH-algebra $H$
automatically has an antipode.  With $f$ its FH-homomorphism
and $t$ a right norm, define $S:  H \rightarrow H$ by 
\begin{equation}
S(a) = \sum f(t_1a)t_2.
\label{eq:antipode}
\end{equation}
Then for every $a \in H$ \[
\sum S(a_1)a_2 = \sum f(t_1 a_1) t_2 a_2 = f(ta) 1 = \epsilon(a) 1.
\] Now in the convolution algebra structure on $\End_k(H)$, this shows
$S$ has $\Id_H$ as right inverse. Since $\End_k(H)$ is
 finite projective over $k$, it follows that 
$\Id_H$ is also a left inverse of $S$; whence $S$ is the unique antipode. 

The Pareigis Theorem \cite{Par71} generalizing the Larson-Sweedler
Theorem \cite{LS} shows that a finite projective
Hopf algebra $H$ over a ground ring $k$ with trivial Picard group
is an FH-algebra.  In detail, the theorem proves the following
in the  order given.  The first two items are proven without the hypothesis
on the Picard group of $k$.  The last two items require only that
$\int^{\ell}_{H^*}$ be free of rank 1. 
\begin{enumerate}
\item There is a right Hopf $H$-module structure on $H^*$.
Since all Hopf modules are trivial, $H^* \cong P(H^*) \otimes H$,
for the coinvariants $P(H^*) = \int^{\ell}_{H^*}$. 
\item The antipode $S$ is bijective.
\item  There exists a left integral $f$ in $H^*$
such that the mapping $\Theta: H \rightarrow H^{\ast}$ defined by
\begin{equation}
\Theta(x)(y) = f(yS(x))
\end{equation}
 is a right Hopf module isomorphism. 
\item $H$ is a Frobenius algebra with 
 Frobenius homomorphism $f$.
\end{enumerate}

It follows from 2. above that an FH-algebra $H$ possesses an
$\epsilon$-invariant anti-automorphism
$S$.  If $f \in H^*$ is an FH-homomorphism, then $Sf$ is a
Frobenius homomorphism and {\it left integral} in $H^*$.  It is therefore
equivalent to replace right with left in Definition~\ref{def-june3}.

Let $m: H \rightarrow k$ be the right modular function
of $H$.  Since $m$ is an algebra homomorphism, it is group-like
in $H^*$, whence $m$ at times is called the {\it right distinguished group-like element}
in $H^*$. 

The next proposition is obtained in an equivalent form
in \cite{OS73}, \cite{FMS} and \cite{BFS}, though in somewhat
different ways. 

\begin{prop} 
Let $H$ be an FH-algebra with FH-homomorphism
 $f$ and right norm 
$t$.  Then $(f, S^{-1}t_2,t_1)$ is 
a Frobenius system for $H$. 
\label{lemma-one}
\end{prop}
\begin{proof} Applying $S^{-1}$ to both sides of Equation~\ref{eq:antipode}
yields 
\begin{equation}
\sum S^{-1}(t_2) f(t_1 a)  =  a, 
\label{eq:june2}
\end{equation}
for every $a \in H$. 
It follows from the finite projectivity assumption on $H$ that
$(f,S^{-1}(t_2),t_1)$ is a Frobenius system.  \end{proof}

It follows from the proposition that $t \leftharpoonup f = 1$.   Together with the corollary below
this implies that $f$ is a right norm in $H^*$, since $1$ is the
counit for $H^*$. It follows  that  $g$ is another 
FH-homomorphism for $H$ iff $g = f\lambda$ for some $\lambda \in k^{\circ}$. 
From Equation~\ref{eq:nesbyen} and the proposition above it follows that
\begin{equation}
S(t) = t \leftharpoonup m.
\end{equation}

\begin{prop}
 $H$ is an FH-algebra if and only if $H^*$ is an FH-algebra. 
\label{cor-june5} 
\end{prop}
\begin{proof}
 It suffices by duality
 to establish the forward implication. Suppose $f$ is an FH-homomorphism
 for $H$ and $t$ a right norm.
Now Equation~\ref{eq:antipode} and the argument after it work
for $H^*$ and the right integrals $t$, $f$ since $t \leftharpoonup f$ is the counit on $H^*$. It follows
that 
\begin{equation}
S(g) = \sum  (f_1 g)(t) f_2
\end{equation}
 is  an equation for the antipode
in $H^*$.  By taking $S^{-1}$ of both sides
we see that \[
(t,S^{-1}f_2,f_1)
\]
 is a Frobenius system
for $H^*$ .  
 Whence $t$ is an FH-homomorphism for $H^*$ with right norm
$f$.  \end{proof}

It follows that $H^*$ is also an augmented Frobenius algebra.
Next, we simplify our criterion for FH-algebra. 
\begin{prop}
If $H$  is an FH-algebra if and only if $H$ is a Frobenius algebra and a Hopf algebra.
\end{prop}
\begin{proof}
The forward direction is obvious.  For the converse, we use the fact that the 
$k$-submodule of integrals of an augmented Frobenius
algebra is free of rank $1$ (cf.\  
\cite[Theorem 3]{Par71} or Proposition~\ref{prop-ein}.  It follows that 
$\int^{\ell}_H \cong k$.  From Pareigis's Theorem we obtain that the dual
 Hopf algebra $H^*$ is a Frobenius algebra. Whence $\int^{\ell}_{H^*} \cong k$
and $H$ is an FH-algebra. 
\end{proof}

Let $b \in H$ be the right {\it distinguished group-like element} satisfying
\begin{equation}
gf = g(b)f
\label{eq:wed}
\end{equation}
 for every $g \in H^*$.  

The convolution product inverse of $m$ is 
$ m^{-1} = m \circ S $. Given a left norm $v \in H$, we claim that   \[
va = vm^{-1}(a).
\]
 Since $t$ is a right norm, 
$S$ an anti-automorphism and $\epsilon$-invariant, 
it follows that $St$ is a left norm.  Then we may assume
$v = St$.  
 Then $S(at) = StSa = m(a) St $, whence
\(
vx = vmS^{-1}(x)
\)
 for every $a,x \in H$. The claim then follows from 
 \(
m \circ S^2 = m, 
\)
since this implies that $m \circ S^{-1} = m^{-1}$. 
But $m \circ S^2 = m^{-1} \circ S = m$, since $m^{-1}$ is group-like.

\begin{lemma}
Given an FH-algebra $H$
with right norm $f \in H^*$ and  right norm $t \in H$ such that $f(t) = 1$, 
the Nakayama automorphism, relative to $f$, and its inverse are
given by:
\begin{equation}
\eta(a) = S^2(a \leftharpoonup m^{-1}) = (S^2a) \leftharpoonup m^{-1}, 
\label{eq:eta(a)}
\end{equation}
\[ \eta^{-1}(a) = S^{-2}(a \leftharpoonup m) = (S^{-2}a) \leftharpoonup m. \]
\label{lemma-june8}
\end{lemma}

\begin{proof} 
Using the Frobenius coordinates $(f,S^{-1}t_2, t_1)$, we note that  
\[
\eta^{-1}(a) = \sum S^{-1}(t_2) f(t_1 \eta^{-1}(a)) = \sum S^{-1}(t_2) 
f(at_1).
\]
We make a computation as in \cite[Lemma 1.5]{FMS}:
\begin{eqnarray*}
S^2(\eta^{-1}(a)) & = & \sum f(at_1) St_2 \\
& = & \sum f(a_1t_1)a_2 t_2 St_3 \\
& = & \sum f(a_1 t) a_2 \\
& = & a \leftharpoonup m 
\end{eqnarray*}
since  $a \leftharpoonup f = f(a)1$,  $at = m(a)t$  for every $a \in H$ 
and $f(t) = 1$.  Whence $\eta^{-1}(a) =
S^{-2}(a \leftharpoonup m)$.  Since $mS^{-2} = m$, it follows
that $\eta^{-1}(a) = (S^{-2}a) \leftharpoonup m$. 

It follows that $a = (S^{-2} \eta a) \leftharpoonup m$, so let the convolution
inverse $m^{-1}$ act on both sides:  $(a \leftharpoonup m^{-1}) = S^{-2} \eta(a)$.
Whence $\eta(a) = S^2(a \leftharpoonup m^{-1}) = (S^2a) \leftharpoonup m^{-1}$,
since $m^{-1} S^2 = m^{-1}$. \end{proof}

As a corollary, we obtain \cite[Folg.\ 3.3]{OS73} and 
\cite[Proposition 3.8]{BFS2}:  
if $H$ is a unimodular
FH-algebra, then the Nakayama automorphism is the square of the
antipode. 

Now recall our definition of $b$ after Proposition~\ref{cor-june5}
as the right distinguished group-like in $H$. 
Equation~\ref{eq:Radford} below was first established in \cite{Rad94}
for finite dimensional Hopf algebras over fields by different means.

\begin{theorem} 
If $H$ is an  FH-algebra with FH-homomorphism
 $f $ and right norm $t$, then 
\begin{equation}
\sum t_2 \otimes t_1 = \sum b^{-1} S^2t_1 \otimes t_2.
\label{eq:Radford}
\end{equation}
\label{prop-gote}
\end{theorem}
\begin{proof}  On the one hand, we have seen that
$(f,S^{-1}t_2,t_1)$ are Frobenius coordinates for $H$.  
On the other hand, the equation $f \rightharpoonup x =  b f(x)$
for every $x \in H$ follows from Equation~\ref{eq:wed} and gives 
\begin{eqnarray*}
\sum S(t_1)b f(t_2 a) &=& \sum  S(t_1) t_2a_1 f(t_3 a_2) \\
& = & \sum a_1 f(t a_2) \\
& = & \sum a_1 \epsilon(a_2) f(t) = a. 
\end{eqnarray*}
Then $(f, S(t_1)b, t_2)$ is another Frobenius system for $H$.

Since $(S^{-1}(t_2), t_1)$ and $(S(t_1)b, t_2)$ are both dual bases to
$f$, it follows  
  that $\sum S^{-1}t_2 \otimes t_1 =  \sum S(t_1)b \otimes t_2$.
 Equation~\ref{eq:Radford} follows from applying $S \otimes 1$
to both sides.  \end{proof}

  Proposition~\ref{lemma-one} with $a = S^{-1}t$
gives \[
\sum S^{-1}t_2 f(t_1 S^{-1}t) = S^{-1} t f(S^{-1}t) = S^{-1}t.
\]
Since $S^{-1}t$ is a left norm, 
it follows that
\begin{equation}
f(S^{-1}t) = 1.
\label{eq:one} 
\end{equation}

The next proposition is not mentioned 
in the literature for
Hopf algebras. 

\begin{prop}
Given an FH-algebra $H$ with
FH-homomorphism $f$,  
the right distinguished group-like element $b$ is equal to the
 derivative $d$ of the left integral Frobenius homomorphism
 $g :=S^{-1}f$  with respect to  $f$:
\begin{equation}
b = \frac{dg}{df}
\end{equation}
\label{prop-shop}
\end{prop}
\begin{proof} By Lemma~\ref{lemma-emma}, another Frobenius system for $H$ is 
given by $(g, St_1,t_2)$, since $S$ is an anti-automorphism. 
Then there exists a (derivative) $d \in H^{\circ}$ such that 
\begin{equation}
df = g.
\label{eq:ena}
\end{equation}

$g$ is a left norm in $H^*$ since $S^{-1}$ is an $\epsilon$-invariant
anti-automorphism.
Also $bf$ is a left integral in $H^*$ by the following argument.
For any $g,g' \in H^*$, we
have $b(gg') = (bg)(bg')$ as $b$ is group-like. Then for every $h \in H^*$
\begin{eqnarray*}
h(bf) & = & b[(b^{-1}h)f] \\
& = & b[(b^{-1}h)(b)f] \\
& = & h(1)(bf). 
\end{eqnarray*}

Now both $g(t)$ and $bf(t)$ equal $1$, since $f(S^{-1}t) = 1$, 
 $f(tb) = \epsilon(b) f(t) = 1$ and $b$ is group-like.
Since $bf$ is a scalar multiple of the norm $g$, it follows that 
\begin{equation}
g = bf.
\label{eq:rad2}
\end{equation}
Finally, $d = b$ since $df = bf$ from Equations~\ref{eq:ena}
and~\ref{eq:rad2}, 
and $f$ is  nondegenerate. \end{proof}

We next give a different derivation for FH-algebras of a  formula in \cite{Rad76}  
for the fourth power of the
antipode of a finite dimensional Hopf algebra.  The main point is
that the Nakayama
automorphisms associated with the two Frobenius homomorphisms 
$S^{-1}f$ and $f$ differ by an inner automorphism determined by the
derivative in Proposition~\ref{prop-shop}. 

\begin{theorem} 
Given an FH-algebra $H$ with right distinguished group-like
elements $m \in H^*$ and $b \in H$, the fourth power of the antipode
is given by 
\begin{equation}
S^4(a) = b(m^{-1} \rightharpoonup a \leftharpoonup m)b^{-1}
\label{eq:friday}
\end{equation}
for every $a \in H$.
\label{th-june8}
\end{theorem}
\begin{proof} Let $g := S^{-1}f$ and denote the left
norm $St$ by  $\Lambda$.  Note that $g(\Lambda) = 1 = g(S^{-1}\Lambda)$
since $f(t) = 1 = f(S^{-1}t)$.
 We note that $(g, \Lambda_2, S^{-1}\Lambda_1)$
are Frobenius coordinates for $H$, since $S$ is an 
anti-automorphism of $H$

Then the Nakayama automorphism $\alpha$ associated with $g$ has inverse
 satisfying
\[
\alpha^{-1}(a)  =  \sum \Lambda_2 g(aS^{-1}\Lambda_1 ) \]
whence 
\begin{eqnarray*}
S^{-1}\alpha^{-1}(a) & = & \sum S^{-1}g(\Lambda_1 Sa) S^{-1}(\Lambda_2) \\
& = & \sum S^{-1}(\Lambda_3)S^{-1}g( \Lambda_1Sa_2)\Lambda_2Sa_1 \\
& = & \sum S^{-1}g(\Lambda Sa_2)Sa_1 \\
& = & g(S^{-1}\Lambda) \sum m^{-1}(Sa_2)Sa_1 = S(m \rightharpoonup a),
\end{eqnarray*}
since $Sm^{-1} = m$.  
It follows that
\begin{eqnarray}
\alpha^{-1}(a) & = S^2(m \rightharpoonup a) & = m \rightharpoonup S^2a \\
\alpha(a) & = m^{-1} \rightharpoonup S^{-2}a & 
= S^{-2}(m^{-1} \rightharpoonup a). 
\end{eqnarray} 

From Proposition~\ref{prop-shop} we have $g = bf = f\eta(b)$,
where $\eta$ is the Nakayama automorphism of $f$.
By  Equation~\ref{eq:alex} and Lemma~\ref{lemma-june8},
\begin{eqnarray*}
m^{-1} \rightharpoonup S^{-2}a &=&  \alpha(a) \\
& = & \eta(b^{-1}) \eta(a) \eta(b) \\
&=& m^{-1}(b^{-1})b^{-1}(S^2(a)\leftharpoonup m^{-1})bm^{-1}(b) \\
& = & b^{-1}S^2(a)b \leftharpoonup m^{-1}, 
\end{eqnarray*}
since $b$ and $m$ are group-likes and $S^2$ leaves $m$ and $b$ fixed.
It follows that 
\[
a = m \rightharpoonup b^{-1}S^4(a)b \leftharpoonup m^{-1}, \]
for every $a \in H$. 
Equation~\ref{eq:friday} follows. \end{proof}

The theorem implies \cite[Corollary 3.9]{BFS2}, which states that
$S^4 = \Id_H$,
  if $H$ and $H^*$ are
unimodular finite projective Hopf algebra over $k$.
For localizing with respect to any maximal ideal ${\mathcal M}$, we obtain unimodular
Hopf-Frobenius
algebras $H_{{\mathcal M}} \cong H \otimes k_{{\mathcal M}}$ and its dual,
since the local ring $k_{{\mathcal M}}$ has trivial Picard group.  
By Theorem~\ref{th-june8},
the localized antipode satisfies $(S_{{\mathcal M}})^4 = \Id$ 
for every maximal ideal
${\mathcal M}$ in $k$; whence $S^4 = \Id_H$
\cite{Sil}. 

\begin{theorem}
Let $k$ be a commutative ring in which $2$ is not a zero divisor
and $H$ a finite projective Hopf algebra.
If $H$ is separable and coseparable, then $S^2 = \Id$.
\label{th-geneg}
\end{theorem}

\begin{proof}
First  we note that $H$ is unimodular and counimodular.
Then it follows from the theorem above that $S^4 = \Id$.
Localizing with respect to the set $T=\{ 2^n ,n=0,1,..\}$ we may assume
that $2$ is invertible in $k$. 
Then $H=H_+ \oplus H_-$ where $H_{\pm} =\{h\in H: S^2 (h)={\pm}h\}$, 
respectively.
We have to prove that $H_- =0$. It suffices to prove that $(H_- )_{\mathcal M} =0$ for
any maximal ideal ${\mathcal M}$ in $k$. Since 
$H_{\mathcal M}/{\mathcal M}H_{\mathcal M}$ is separable and coseparable
over the field $k/{\mathcal M}$, we deduce from the main theorem in 
\cite{EG} that  $(H_- )_{\mathcal M} \subset {\mathcal M}H_{\mathcal M}$ and therefore  
$(H_-)_{\mathcal M} \subset {\mathcal M}(H_- )_{\mathcal M} $. 
The desired result follows from the
Nakayama Lemma because $H_-$ is a direct summand in $H$.
\end{proof}

In \cite{BFS2} it was established that if $H$ is separable over a ring
$k$ with no torsion elements, then $S^2 = \Id$. We may improve on this
and similar results by an application of  the last theorem.
If $k$ is a commutative ring
and $M$ is a finitely generated projective $k$-module, we let $rank_M:\, 
Spec\, k \rightarrow {\mathcal Z }$ be the \textit{rank function},
which
is defined  on a prime ideal ${\mathcal P}$ in $k$ by 
$$rank_M ({\mathcal P}):=
dim_{\ \overline{k/{\mathcal P}}}\ (M\otimes_k k/{\mathcal P})\otimes_{k/{\mathcal P}}
{\overline{k/{\mathcal P}}}$$
where ${\overline{k/{\mathcal P}}}$ is the field of fractions of ${k/{\mathcal
P}}$.
The range of $rank_M$ is finite and consists of a set of positive
integers $n_1 , n_2, ... ,n_k$.

Now for any prime $p \in \mathcal{Z}$, let $Spec^{(p)}k \subseteq Spec\, k$
be the subset of prime ideals $P$ for which the characteristic
$char\, \overline{(k/P)} = p$.  Suppose that $Spec^{(p)}k$ is non-empty
and $$ rank_M(Spec^{(p)}k) = 
\{ n_{i_1},\ldots, n_{i_s} \}. $$
For such $p$ and $\phi$ the Euler function, we define 
$$ N(M,p) := \max_{m = 1,\ldots,s} \{ n_{i_m}^{\frac{\phi(n_{i_m})}{2}} \, \} .$$
  
\begin{theorem}
Let $k$ be a commutative ring in which $2$ is not a zero divisor
and $H$ be a f.g. projective Hopf algebra.  If $H$ is $k$-separable such that
$ N(H,p) < p$ for every odd prime $p$, then $H$ is coseparable and $S^2 = \Id$.
\label{th-saturn}
\end{theorem}
\begin{proof}
First we note that $2$ may be assumed invertible in $k$ without loss of
generality by localization with respect to powers of $2$.
Let $\mathcal{M}$ in $k$ be  a maximal ideal. The characteristic
of $k/\mathcal{M}$ is not $2$ by our assumption. 
 
It is known that an algebra $A$ is separable iff $A/\mathcal{M}A$
is separable over $k/\mathcal{M}$ for every maximal ideal $\mathcal{M}
\subset k$ \cite{DM}: whence $H/\mathcal{M}H$ is $k/\mathcal{M}$-separable. 
Furthermore note that if $d(\mathcal{M}) := \dim_{k/\mathcal{M}}
H/\mathcal{M}H $ is greater than $2$, then 
$$
d(\mathcal{M})^{\frac{\phi(d(\mathcal{M}))}{2}} < char\, (k/\mathcal{M}). $$
It then follows from \cite{EG}
that $H^*/\mathcal{M}H^* \cong (H/\mathcal{M}H)^*$ is $k/\mathcal{M}$-separable
for such $\mathcal{M}$.

If $d(\mathcal{M}) = 2$ and $\overline{k/\mathcal{M}}$ denotes the algebraic closure
of $k/\mathcal{M}$, then $H^*/mH^* \otimes_k \overline{k/\mathcal{M}}$ is either
semisimple or isomorphic to the ring of dual numbers.  But the latter is impossible
since it is not a Hopf algebra in characteristic different from $2$.  
Hence $H^*/\mathcal{M}H^*$ is $k/\mathcal{M}$-separable
for all maximal ideal $\mathcal{M}$.  Hence $H^*$ is $k$-separable
by \cite{DM}.  
By Theorem~\ref{th-geneg} then,
$S^2 = \Id$. 
\end{proof}

In closing this section, we note that 
Schneider \cite{Arg} has  established 
Equation~\ref{eq:friday}  by different 
methods for $k$ a field. 
Equation~\ref{eq:friday} is generalized
 in a different
direction by Koppinen \cite{Kop98}. Waterhouse  sketches a different method of how to 
extend the Radford formula to a 
finite projective Hopf algebra \cite{Water}.  

\section{FH-Subalgebras}

In this section we prove 
that a Hopf subalgebra pair of FH-algebras $B \subseteq A$ 
form a  $\beta$-Frobenius extension. The first results of this kind were obtained
by Oberst and Schneider in \cite{OS73} under the assumption that $H$ is cocommutative. 

The proposition below {\it sans} Equation~\ref{eq:99}
 is more general than  \cite[Theorem 1.3]{Far} and a special case of \cite[Satz 7]{Par64} : the
proof simplifies somewhat and is needed for establishing Equation~\ref{eq:99}. 

\begin{prop}
Suppose $A$ and $B$ are Frobenius algebras over the same commutative
ring $k$ with Frobenius coordinates $(\phi,x_i,y_i)$
and $(\psi,z_j,w_j)$, respectively. If $B$ is a subalgebra of $A$ such
that  $A_B$ is projective and the Nakayama automorphism
$\eta_A$ of $A$ satisfies $\eta_A(B) = B$,
then $A/B$ is a $\beta$-Frobenius extension with $\beta$ the relative
Nakayama automorphism, 
\begin{equation}
\beta =  \eta_B \circ \eta_A^{-1},
\end{equation}
and $\beta$-Frobenius homomorphism
\begin{equation}
F(a) = \sum_j \phi(az_j)w_j, 
\label{eq:99}
\end{equation}
for every $a \in A$. 
\label{prop-seven}
\end{prop}
\begin{proof} Since $B$ is finite projective over $k$, it follows that
$A_B$ is a {\it finite} projective module. 

     It remains to check that ${}_BA_A \cong {}_{\beta}(A_B)^*_A$,
which we do below by using the Hom-Tensor Relation and Equation~\ref{eq:yama}
twice (for $A$ and for $B$). Let $\eta_A^{-1}$ denote the restriction 
of $\eta_A^{-1}$ to $B$ below. 
\begin{eqnarray*}
{}_BA_A & \cong & {}_{\eta^{-1}_A} \Hom_k( A,k)_A \\
& \cong & \Hom_k(A \otimes_B B_{\eta^{-1}_A}, k)_A \\
& \cong & \Hom_B(A_B, {}_{\eta^{-1}_A}\Hom_k(B,k)_B)_A \\
& \cong & {}_{\eta^{-1}_A}\Hom_B(A_B, {}_{\eta_B}B_B) \\
& \cong & {}_{\eta_B \circ \eta^{-1}_A}\Hom_B(A_B,B_B)_A .  
\end{eqnarray*}

By sending $1_A$ along the isomorphisms in the last set of equations,
we compute that the Frobenius homomorphism $F: {}_BA_B \rightarrow
{}_{\beta}B_B$ is given by Equation~\ref{eq:99}. 
One may double check that $F(bab') = \beta(b)F(a)b'$ for every $b,b' \in B,
a \in A$ by applying Equation~\ref{eq:laurdag}.
\end{proof}

Given a commutative ground ring $k$, we assume $H$ and $K$  
are Hopf algebras with $H$ a finite projective $k$-module. 
$K$ is a \textit{Hopf subalgebra}  of $H$ if it is
a \textit{pure} $k$-submodule of $H$ \cite{Lam} and a subalgebra of $H$ 
for which $\Delta(K) \subseteq K \otimes_k K$ and 
$S(K) \subseteq K$. It follows that $K$ is 
finite projective as a $k$-module \cite{Lam}. The next lemma is a corollary
of the Nicholls-Zoeller freeness theorem. 

\begin{lemma}
If $H$ is a finitely generated  free Hopf algebra over a local ring $k$
with $K$ a Hopf subalgebra, then the natural modules $H_K$ and ${}_KH$ are free.
\label{lemma-topfjell}
\end{lemma}
\begin{proof}
It will suffice to prove that $H_K$ is free, the rest of the proof being entirely similar. 
First note that $H_K$ is finitely generated since $H_k$ is. 
If $ \mathcal{M}$ is  the maximal ideal of $k$, then the finite dimensional
Hopf algebra $\overline{H} := H/\mathcal{M}H$ is  free over the Hopf subalgebra 
$\overline{K} := K/\mathcal{M}K$
by purity and the freeness theorem in \cite{NZ}. Suppose $\theta: \overline{K}^n \stackrel{\cong}{\rightarrow}
\overline{H}$ is a $\overline{K}$-linear isomorphism.  Since $K$
is finitely generated over $k$, 
$\mathcal{M}K$ is contained in the radical of $K$. 
Now  $\theta$ lifts
to a right $K$-homomorphism $K^n \rightarrow H$ with respect to the natural projections
$H \rightarrow \overline{H}$ and $K^n \rightarrow \overline{K}^n$. 
By Nakayama's lemma, the homomorphism $K^n \rightarrow H$ is  epi (cf.\ \cite{Sil}).  
Since $H_k$ is finite projective, $\tau$ is a $k$-split epi, 
which is bijective
by Nakayama's lemma  applied to the underlying $k$-modules. 
Hence, $H_K$ is free of finite rank. 
\end{proof}

Over a non-connected ring $k = k_1 \times k_2$,
it is easy to construct examples of   Hopf subalgebra pairs
\[
K := k[H_1 \times H_2] \subseteq H := k[G_1 \times G_2]
\]
 where $G_1 > H_1$,
$G_2 > H_2$ are  subgroup pairs of finite groups and $H_K$ is not free 
(by counting dimensions on either side of $H \cong K^n$). 
  The next proposition follows from the lemma. 

\begin{prop}
If $H$ is a finite projective Hopf algebra and $K$ is a 
finite projective Hopf subalgebra of $H$, then the natural modules $H_K$ and ${}_KH$
are finite projective. 
\label{lemma-fjelltop}
\end{prop}
\begin{proof}
 We prove only that $H_K$ is finite projective since the proof that ${}_KH$ 
is entirely similar. 
 First note that $H_K$ is finitely generated.

If $k$ is a commutative ground ring, $Q \rightarrow P$ is an epimorphism
of $K$-modules, then it will suffice to show
that the induced map $\Psi: \Hom_K(H_K, Q_K) \rightarrow \Hom_K(H_K,P_K)$ is
epi too. Localizing at a maximal ideal
$\mathcal{M}$ in $k$, we obtain a homomorphism denoted by $\Psi_{\mathcal M}$. 
By adapting a standard argument such as in \cite{Sil}, we note that for every module $M_K$
\begin{equation}
\Hom_K(H_K,M_K)_{\mathcal{M}} \cong \Hom_{K_{\mathcal{M}}}^{r}(H_{\mathcal{M}},
M_{\mathcal{M}})
\label{eq:isom}
\end{equation}
since $H_k$ is finite projective.  Then $\Psi_{\mathcal M}$
maps $$ \Hom_{K_{\mathcal{M}}}^{r}(H_{\mathcal{M}},
Q_{\mathcal{M}}) \rightarrow \Hom_{K_{\mathcal{M}}}^{r}(H_{\mathcal{M}},
P_{\mathcal{M}}).$$

By Lemma~\ref{lemma-topfjell}, $H_{\mathcal{M}}$ is free over
$K_{\mathcal{M}}$.  It follows that $\Psi_{\mathcal{M}}$ is epi
for each maximal ideal $\mathcal{M}$, whence $\Psi$ is epi. 
\end{proof}

Suppose $K \subseteq H$ is a  pair of FH-algebras where $K$ is a Hopf
subalgebra of $H$:  call $K \subseteq H$ a 
{\it FH-subalgebra pair}.  We now easily prove
 that $H/K$ is a $\beta$-Frobenius extension.  

\begin{theorem}
If $H/K$ is a FH-subalgebra pair, then $H/K$ is a $\beta$-Frobenius extension 
where $\beta = \eta_K \circ \eta_H^{-1}$.
\label{theorem-grey}
\end{theorem} 
\begin{proof}
The Nakayama automorphism $\eta_H$ sends $K$ into $K$ by 
Equation~\ref{eq:eta(a)}, since $K$ is a Hopf subalgebra of $H$. 
 $H_K$ is projective by  Proposition~\ref{lemma-fjelltop}. 
The conclusion follows then from Proposition~\ref{prop-seven}
\end{proof}

From the theorem and Lemma~\ref{lemma-june8} we readily
compute $\beta$ in terms of the relative
modular function $\chi := m_H * m_K^{-1}$,
obtaining the formula \cite[1.6]{FMS}: for every $x \in K$, 
\begin{equation}
\beta(x) = x \leftharpoonup \chi.
\label{eq:snow}
\end{equation}

Applying $m_K$ to both sides of this equation, we obtain
\begin{equation}
m_H(x) = m_K(\beta(x)), 
\label{eq:improved}
\end{equation}
 a formula which extends that in \cite[Corollary 1.8]{FMS} 
from the case $\beta = \Id_K$.

\section{Some Formulas for a Hopf Subalgebra Pair} 

It follows
from Theorem~\ref{theorem-grey} and   Lemma~\ref{lemma-topfjell}
 that a Hopf subalgebra pair  $K \subseteq H$  {\it over a local ring} $k$ is a free 
$\beta$-Frobenius extension.  Since $H_K$ is free and therefore faithfully flat,
the proof in \cite{FMS} that $(E, S^{-1}(\Lambda_2),\Lambda_1)$, 
defined below, 
is a Frobenius system carries through word for word as described next.

From Proposition~\ref{lemma-one} it follows
that $(f, S^{-1}(t_{H(2)}),t_{H(1)})$ is a  Frobenius system for $H$
where $f \in H^*$ and  $t_H$ in $H$ are right integrals
 such that $f(t_H) = 1$.
Given  right and left modular functions $m_H$ and $m_H^{-1}$,
a  computation using Equation~\ref{eq:formue} 
determines that   
\begin{equation}
\eta_H^{-1}(a) = S^{-2}(a \leftharpoonup m_H),
\label{eq:soon}
\end{equation}
for every $a \in H$. 
Let $t_K$  be a right integral for $K$.  Now by a theorem  in \cite{NZ}, 
 $H_K$ and
${}_KH$ are free. Then there exists $\hat{\Lambda} \in H$ such that
$t_H = \hat{\Lambda}t_K$.  Let $\Lambda := \eta_H(S^{-1}(\hat{\Lambda}))$.
Then a $\beta$-Frobenius system for $H/K$ is given
by $( E, S^{-1}\Lambda_{(2)}, \Lambda_{(1)} )$ where
\begin{equation}
E(a) = \sum_{(a)} f(a_{(1)} S^{-1}(t_K))a_{(2)}, 
\label{eq:fishmont}
\end{equation}
for every $a \in H$ \cite{FMS}.  For example,
if $K$ is the unit subalgebra, $E = f$ and $\Lambda = t$. 

The rest of this section is devoted to comparing the different Frobenius systems
for  a Hopf subalgebra pair $K \subseteq H$  over a local ring $k$
implied by our work in Sections~4 and~5. 
Suppose that $f \in \int^r_{H^*}$ and $t \in \int^r_H$ such that $ft = \epsilon$,
and that $g\in \int^r_{K^*}$ and $n \in \int^r_K$ satisfy $gn = \epsilon |_K$.
Then by Section 4 $(f, S^{-1}(t_2),t_1)$ is a Frobenius system for $H$,
and $(g,S^{-1}(n_2),n_1)$ is a Frobenius system for $K$, both as Frobenius
algebras. 

By Equation~\ref{eq:99}, we note that a Frobenius
homomorphism $F: H \rightarrow K$ of the $\beta$-Frobenius
extension $H/K$ is  given by
\begin{equation}
F(a) =  \sum f(aS^{-1}(n_2))n_1.
\end{equation}
Comparing $E$ in Equation~\ref{eq:fishmont} and
$F$ above, we  compute  the (right) derivative
$d$ such that $F = Ed$: 
\begin{eqnarray*}
d & = & \sum F(S^{-1}(\Lambda_2)\Lambda_1 \\
& = & \sum f(S^{-1}(\Lambda_2) S^{-1}(n_2)) n_1 \Lambda_1 \\
& = & \sum (S^{-1}f)(n_2\Lambda_2)n_1 \Lambda_1 = (S^{-1}f)(n\Lambda)1_H, 
\end{eqnarray*}
since $S^{-1}f \in \int^{\ell}_{H^*}$. Hence, 
$(S^{-1}f)(n\Lambda) \in k^{\circ}$.  

We  next make note of a transitivity lemma for Frobenius systems, 
which adds Frobenius systems to the transitivity theorem, \cite[Satz 6]{Par64}.

\begin{lemma}
Suppose $A/S$ is a $\beta$-Frobenius extension with system $(E_S,x_i,y_i)$
and $S/T$ is a $\gamma$-Frobenius extension with system $(E_T,z_j,w_j)$.
If $\beta(T) = T$, then $A/T$ is a $\gamma \circ \beta$-Frobenius extension
with system,
\[
(E_T \circ E_,\ x_i z_j,\ \beta^{-1}(w_j)y_i).
\]
\label{lemma-femme}
\end{lemma}
\begin{proof}
The mapping $E_TE_S$ is clearly a bimodule homomorphism ${}_T A_T \rightarrow {}_{\gamma \circ \beta}T_T$.
We compute for every $a \in A$:
\begin{eqnarray*}
\sum_{i,j} x_iz_j E_TE_S(\beta^{-1}(w_j)y_ia) & = & \sum_i x_i \sum_j z_j E_T(w_j E_S(y_ia)) \\
& = & \sum_i x_i E_S(y_ia) = a, \\
\sum_{i,j} (\gamma \beta)^{-1}(E_TE_S(ax_iz_j)) \beta^{-1}(w_j)y_i & = & 
\sum_i \beta^{-1}(\sum_j \gamma^{-1}(E_T(E_S(ax_i)z_j))w_j) y_i \\
& = & \sum_i \beta^{-1}(E_S(ax_i))y_i = a. \qed
\end{eqnarray*}
\renewcommand{\qed}{}\end{proof}

Applying the lemma to the Frobenius system $(E,S^{-1}(\Lambda_2),\Lambda_1)$ for 
$H/K$ and Frobenius system $(g,S^{-1}(n_2),n_1)$ for
$K$ yields the  Frobenius system for the algebra $H$, 
\[
( g \circ E,\ S^{-1}(\Lambda_2)S^{-1}(n_2),\ \beta^{-1}(n_1) \Lambda_1).\]
Comparing this with the Frobenius system $(f,S^{-1}(t_2),t_1)$,
we compute the derivative $d' \in H^{\circ}$ such that
$(gE)d' = f$:
\begin{equation}
d' = \sum f(S^{-1}(n_2 \Lambda_2)) \beta^{-1}(n_1)\Lambda_1
\end{equation}
We note that $f = gF$, since for every $a \in H$, 
\[
g(\sum f(aS^{-1}(n_2))n_1) = f(a \sum S^{-1}(n_2)g(n_1)) =  f(a).
\]
Now appy $g$ from the left to $F = Ed$ and conclude that  
 $d = d'$.  It follows  that $gE$ is a right norm in $H^*$ like $f$,
since $d \in k^{\circ}$. 

Since $ft = \epsilon$ and $m_H(d') = d  = (S^{-1}f)(n\Lambda)1_H$, 
we see that $dt$ is a right norm for $gE$.
Using Equation~\ref{eq:maritshoes}, we compute that
\begin{eqnarray}
dt & = & \sum \epsilon(S^{-1}(n_2\Lambda_2))\beta^{-1}(n_1)\Lambda_1  \nonumber \\
& = & \beta^{-1}(n) \Lambda. 
\label{eq:furelise}
\end{eqnarray}
Recalling from Section 1 that $t = \hat{\Lambda} n$,
we note that 
\begin{equation}
\beta^{-1}(n) \Lambda = \hat{\Lambda} n d.
\label{eq:cecilia}
\end{equation} 

Multiplying both sides of the equation $\beta^{-1}(n) \Lambda = t d$
from the left by $\beta^{-1}(x)$, where $x \in K$, derives Equation~\ref{eq:improved}
by other means for local ground rings.

\section{Finite order elements}

Let $M$ be a finite projective module over a commutative ring 
$k$. Let $rank_M :Spec(k)\to {\mathcal Z}$ be the  rank function as
in Section~4.  
We introduce
the {\it rank number} $\hat{D}(M,k)$ of $M$ as the least common multiple of 
the integers in the range of the rank function on $M$:
\[
\hat{D}(M,k) = l.c.m.\{n_1 , n_2,\ldots ,n_k \}.
\]

Let $H$ be a finite projective Hopf algebra over a Noetherian ring 
$k$.
Let $d \in H$ be a group-like
element.  
 In this section we provide a proof that 
$d^N = 1$
where
$N$ divides
$\hat{D}(H,k)$ (Theorem~\ref{th-GH}). 
In particular if $H$ has constant rank $n$,  such as when 
$Spec(k)$ is connected, then $N$ divides $n$. 
Then we  establish in Corollaries~\ref{cor-I} and~\ref{cor-II}
 that the antipode
 $S$ and Nakayama automorphism $\eta$ satisfy
$S^{4N}={\eta}^{2N}= \Id_H$ as corollaries of Theorem~\ref{th-june8}. 

 Let $k[d,d^{-1}]$
denote the subalgebra of $H$ generated over $k$
 by $1$ and the negative and positive powers of $d$.
Let $k[d]$ denote only the $k$-span of $1$ and the positive powers of
$d$.  Clearly $k[d,d^{-1}]$ is 
 Hopf subalgebra of 
  $H$.  $d $ has a {\it minimal polynomial}
 $p(x) \in k[x]$  if $p(x)$ is a 
  polynomial of least degree such that $p(d) = 0$ and the
$\gcd$ of all the coefficients is $1$. We first consider the case
where $k$ is a domain.

\begin{lemma}
If $k$ is a domain,  each group-like
$d \in H$ has a minimal polynomial of the form $p(x) = x^s - 1$
for some integer $s$. 
Moreover, $s$ divides $dim_{\overline{k}}
(H\otimes_k {\overline{k}})$
and ${\overline f}(d^s)\neq 0$, where
$\overline{k}$ denote the field of fractions of $k$
${\overline f}$ is FH-homomorphism for ${\overline k}[d,d^{-1}]$.
\label{lemma-samedi}
\end{lemma}

\begin{proof} We work at first in the Hopf algebra 
$H \otimes_k \overline{k}$ in which
$H$ is embedded.  Since $\overline{k}[d,d^{-1}]$ is a finite dimensional
Hopf algebra,  there is a unique minimal polynomial of $d$,
 given by 
  $\overline{p}(x) = x^s + \lambda_{s-1} x^{s-1} + \cdots + \lambda_0 1$.  Since
$d$ is invertible, $\lambda_0 \neq 0$ and $\overline{k}[d,d^{-1}] = 
\overline{k}[d]$.  

$\overline{k}[d]$ is a Hopf-Frobenius
algebra with FH-homomorphism $f: \overline{k}[d] 
\rightarrow
\overline{k}$.  Then $f(d^k) d^k = f(d^k) 1$ for every integer $k$,
since each $d^k$ is group-like.  If $f(d^k) \neq 0$, then
$k \geq s$, since otherwise $d$ is root of
 $x^k - 1$, a polynomial of degree less than $s$.
   
Thus, $f(d) = \cdots = f(d^{s-1}) = 0$, but $f(1) \neq 0$ since
$f \neq 0$ on $\overline{k}[d]$. Then $f(\overline{p}(d)) = f(d^s) + \lambda_0 f(1)= 0$,
so that $f(d^s) = -\lambda_0f(1) \neq 0$.
Since $f(d^s)d^s = f(d^s) 1$, it follows that $d^s - 1 = 0$. 
Clearly $\overline{k}[d]$ is a Hopf subalgebra of $H \otimes_k
\overline{k}$ of dimension $s$ over $\overline{k}$
and it follows from the Nichols-Zoeller theorem that $s$
divides  $dim_{\overline{k}}(H\otimes_k {\overline{k}})$.

For $H$ over an integral domain we arrive instead
at  $r(d^s - 1) = 0$ for
some $0 \neq r\in k$.
Since $H$ is finite projective over an integral domain,
it follows that $d^s - 1 = 0$. 
\end{proof}

It follows easily from the proof that if $g(x) \in k[x]$ such that $g(d) = 0$,
then $d^s = 1$ for some integer $s \leq \deg g$. 

\begin{theorem}
Let $H$ be a finite projective Hopf algebra over a commutative ring
$k$, which  contains no additive torsion elements. 
If $d \in H$ is a group-like element, then $d^N = 1$ 
for some $N$ that divides $\hat{D}(H,k)$.
\label{th-samedi}
\end{theorem}

\begin{proof}
  Let ${\mathcal P}$ be a prime ideal in $k$
and $rank_H ({\mathcal P})=n_i$.  Let $D = \hat{D}(H,k)$.

Note that 
$H/{\mathcal P}H \cong H \otimes_k (k/{\mathcal P})$ is a finite 
projective Hopf algebra over the domain $k/{\mathcal P}$.  
By Lemma~\ref{lemma-samedi}, there is an integer $s_{{\mathcal P}} $
such that \(
d^{s_{{\mathcal P}}} - 1 \in {\mathcal P}H
\) 
and $s_{{\mathcal P}} $
divides $n_i$.
It follows that \[
d^{D} - 1 \in {\mathcal P}H \]
for each prime ideal ${\mathcal P}$ of $k$. 
Since $H$ is a finite projective over $k$, a standard argument
gives ${\rm Nil}(k) H = \cap({\mathcal P}H)$ over all prime ideals, 
 where the nilradical ${\rm Nil}(k) = \cap {\mathcal P}$ is equal to the intersection
of all prime ideals in $k$. Thus, $d^{D} - 1 = \sum r_i a_i$ where
$r_i \in {\rm Nil}(k)$.  Let $k_i$ be integers such that $r_i^{k_i} = 0$.
Then 
\begin{equation}
(d^{D} - 1)^{(\sum_{i=1}^n k_i) + 1} = 0.
\label{eq:alexis}
\end{equation}
It is clear that $P(x) :=  (x^{D} - 1)^{(\sum_{i=1}^n k_i) + 1}$ is a
monic polynomial with integer coefficients.

In general, let $m\in {\mathcal Z}$ 
be the least number such that $m\cdot 1=0$: we
 will call $m$ the characteristic of $k$. Clearly in this case
$m=0$ and
${\mathcal Z} \subseteq k$.
Again by Equation~\ref{eq:alexis}, ${\mathcal Z}[d,d^{-1}] = {\mathcal Z}[d]$
is a Hopf algebra over ${\mathcal Z}$. Moreover, since $k$ has no additive
torsion elements, we conclude
that ${\mathcal Z}[d]$ is a free module over ${\mathcal Z}$.

Now by Lemma~\ref{lemma-samedi} $q(x)=x^s - 1$ 
is the minimal polynomial for $d$ over ${\mathcal Z}$ and therefore 
$x^s - 1$ divides $P(x)$ in ${\mathcal Z}[x]$. Since $P(x)$ has the same
roots in ${\mathcal C}$ as $x^{D} - 1$ it follows that a primitive
$s$-root of unity is a $D$-root of unity and whence 
$s$ divides  $D$.
\end{proof}


In preparation for the next theorem, observe that if $(k, {\mathcal M})$ is a local
ring of positive characteristic, then ${\it char}(k)=p^t$  
for some positive power of a prime number $p$  and ${\it char}(k/{\mathcal M})=p$.

\begin{theorem}
Let $k$ be a local ring of positive characteristic and
$H$ be a finite projective Hopf algebra over $k$ of rank 
$n$. If $d\in H$ is group-like, then $d^s=1$ where   $s$
is the order of the image of $d$ in $H/{\mathcal M}H$ (and therefore
$s$ divides $n$). 
\label{th-loc}
\end{theorem}
\begin{proof}
 ${\mathcal Z}_{p^t}$ is clearly a  subring of $k$.
Since we can choose a spanning set of $H$ over $k$ of the 
form $1,d, \ldots,  d^{s-1},t_s,\cdots t_n$
where the elements of the set are linearly independent
{\it modulo} ${\mathcal M}$, it follows that
$1,d,\ldots, d^{s-1}$ generate a free module over  ${\mathcal Z}_{p^t}$.
We claim that  ${\mathcal Z}_{p^t}[d]$ coincides with this module
and therefore is free over ${\mathcal Z}_{p^t}$.

First we  need to prove that ${\mathcal M}H\cap{\mathcal Z}_{p^t}[d]=
p{\mathcal Z}_{p^t}[d]$. To do this, we observe that 
${\mathcal M}\cap {\mathcal Z}_{p^t}=p{\mathcal Z}_{p^t}$ and 
thus $p{\mathcal Z}_{p^t}[d]\subset {\mathcal M}H\cap{\mathcal Z}_{p^t}[d]$.
Then there is a canonical epimomorphism of
algebras over ${\mathcal Z}_{p}:$ ${\frac{{\mathcal Z}_{p^t}[d]}{p{\mathcal Z}_{p^t}}}
\to {\frac{{\mathcal Z}_{p^t}[d]}{ {\mathcal M}H\cap{\mathcal Z}_{p^t}[d]}}$.
Let ${\bar d}$ denote the image of $d$ in $H/{\mathcal M}H$.
Since ${\bar d}^s=1$ over $k/{\mathcal M}$, we  deduce
that 
${\bar d}^s=1$ over ${\mathcal Z}_p$ from the fact that  ${\mathcal Z}_p$ is the prime subfield
of $k/{\mathcal M}$. 
Thus there is a canonical algebra epimorphism
${\mathcal Z}_p [\pi_s]\to {\frac{{\mathcal Z}_{p^t}[d]}{p{\mathcal Z}_{p^t}}}$,
where $\pi_s$ is a cyclic group of order $s$. Since $1,d, \ldots, 
d^{s-1}$
are linearly independent {\it modulo} ${\mathcal M}$,  it follows that
$dim_{{\mathcal Z}_p} ({\frac{{\mathcal Z}_{p^t}[d]}{ {\mathcal M}H\cap{\mathcal
Z}_{p^t}[d]}})\geq s$ while $dim_{{\mathcal Z}_p} ({\mathcal Z}_p [\pi_s])=s$.
Therefore all three ${\mathcal Z}_p$-algebras above have dimension
$s$ and ${\mathcal M}H\cap{\mathcal Z}_{p^t}[d]=
p{\mathcal Z}_{p^t}[d]$.

The next step we need is to prove that $d$ satisfies a monic polynomial
equation of degree $s$ over $k$. Clearly $d^s-1\in 
{\mathcal M}H\cap{\mathcal Z}_{p^t}[d]$
and hence $d^s-1=\Sigma a_i d^i$ with $a_i\in p{\mathcal Z}_{p^t}$.
If all $i<s$ then this is exactly what we need. Otherwise we can
replace $a_{s+k} d^{s+k}$ by $a_{s+k}d^k +a_{s+k}\Sigma a_i d^{i+k}$.
However new coefficients $a_{s+k}a_i$ are divisible by $p^2$.
Continuing this process we will arrive at a monic polynomial in $d$
of degree $s$ because $p^t=0$ and all the monomials 
of degree greater than $s$ will be eliminated.

Now it follows that ${\mathcal Z}_{p^t}[d,d^{-1}]={\mathcal Z}_{p^t}[d]$
is a Hopf algebra over ${\mathcal Z}_{p^t}$ and a free module of rank
$s$. Then ${\mathcal Z}_{p^t}[d]$ is a Hopf-Frobenius algebra because
${\mathcal Z}_{p^t}$ is a local ring. Let $f$ be a Hopf-Frobenius 
homomorphism for ${\mathcal Z}_{p^t}[d]$. If ${\bar f}=f {\rm mod}\ p$,
it is clear that this 
is a Hopf-Frobenius homomorphism for ${\mathcal Z}_{p}[{\bar d}]$. Since
${\bar f}({\bar d}^s)\neq 0$, it follows that
$f(d^s)$ is an invertible element of ${\mathcal Z}_{p^t}$. Hence, 
the relation $f(d^s)d^s=f(d^s)\cdot 1$ implies that $d^s=1$, which 
proves the theorem. 
\end{proof}

\begin{definition}
We say that a commutative ring $k$ is a 
 GH-ring if any group-like
element $d$ of any finite projective Hopf algebra $H$ satisfies
$d^{\hat{D}(H,k)}=1$.
\end{definition}

We have already proved that fields, rings without additive torsion,
and local rings of positive characteristic are GH-rings. Next we make
the following easy remarks:
\begin{remark}
If $k_1,\ldots,k_n$ are GH-rings then $\oplus_{i=1}^n k_i$ is
a GH-ring.
\end{remark}
Let $g:k\to K$ be a ring homomorphism 
($g(1)=1$) and let $g^*:Spec(K)\to
Spec(k)$ be the induced continuous mapping. Then it is well-known
(see for instance \cite{Ba}) that 
$rank_{M\otimes_k K}=g^*\circ rank_M$ for any projective 
$k$-module $M$
and it follows that $\hat{D}(M\otimes_k K,K)$ divides $\hat{D}(M,k)$. 
Then we can make the following
\begin{remark}
If $g:k\to K$ is an embedding and $K$ is a GH-ring, then $k$ is
a GH-ring.
\end{remark}
\begin{theorem}
Noetherian rings are GH-rings.
\label{th-GH}
\end{theorem}
\begin{proof}
 Let $k$ be a Noetherian ring and $T(k)\subset k$ be the set
of all torsion elements, i.e. for any $a\in T(k)$ there exists
a positive integer $m$ such that $ma=0$.  $T(k)$ is clearly an ideal
of $k$. Since $T(k)$ is finitely generated over $k$, 
there exists a positive integer $t(k)$ such that $t(k)T(k)=0$.
Let $\pi_1 :k\to {\frac{k}{T(k)}}$ and $\pi_2 :k\to 
{\frac{k}{t(k)\cdot k}}$ be
canonical surjections.
We claim that $\pi_1 \oplus \pi_2 :k\to  
{\frac{k}{T(k)}}\oplus {\frac{k}{t(k)\cdot k}}$
is an embedding. Indeed, if $(\pi_1 \oplus \pi_2)(x)=0$ then
$x\in T(k)$ and $x=t(k)a$ for some $a\in k$. Obviously then $a\in
T(k)$
and $x=t(k)a=0$. 

Since ${\frac{k}{T(k)}}$ has no additive torsion and therefore is a
GH-ring, it remains
to prove that a Noetherian ring of a positive characteristic is
a GH-ring. For that let us consider a multiplicatively closed set
${\mathcal S}$ consisting of all the  non-divisors of zero of $k$.
It is well-known that $k\to {\mathcal S}^{-1}k$ is an embedding and
$ {\mathcal S}^{-1}k$ is a semi-local ring if $k$ is Noetherian (see
\cite{Ba}). 
So, it is sufficient to prove that a semi-local ring $A$ of positive
characterestic is a GH-ring. Let ${\mathcal M}_1,.., {\mathcal M}_n$ be the
set of all maximal ideals of $A$ and $A_{{\mathcal M}_i}$ be the
corresponding
localizations. Now let us consider a homomorphism 
$f: A\to \oplus A_{{\mathcal M}_i}$
induced by canonical homomorphisms $f_i :A\to A_{{\mathcal M}_i}$.
We claim that $f$ is an embedding. Let in contrary $f(x)=0$.
Then $f_i (x)=0$ for any $i$ and there exists 
$a_i\in A\backslash {{\mathcal M}_i}$
such that $a_i x=0$. Let us consider the ideal $I$ generated by
all $a_i$. Clearly $Ix=0$. On the other hand $I$ cannot belong to
${\mathcal M}_i$ because $a_i$ is not in ${\mathcal M}_i$. Therefore 
we get that $I=A$ and consequently $x=0$. Since any  
$A_{{\mathcal M}_i}$ is a local ring of positive characteristic,
Theorem 3.2 implies the required result.
\end{proof}

As a consequence of Proposition~\ref{th-june8}, Theorem~\ref{th-samedi}
and Equation~\ref{eq:eta(a)}, 
we obtain the following corollaries.

\begin{corollary}
Let $H$ be an FH-algebra
over a Noetherian ring $k$.  Then $S^{4\hat{D}(H,k)} = \eta^{2\hat{D}(H,k)} = \Id_H$
\label{cor-I}
\end{corollary}
\begin{proof}
 Note that $\hat{D}(H,k)=\hat{D}(H^*,k)$.  
\end{proof}

\begin{corollary}
Let $H$ be a finite projective Hopf algebra over a Noetherian ring $k$.
Then $S^{4\hat{D}(H,k)} = \Id_H$.
\label{cor-II}
\end{corollary}
\begin{proof}
 Localizing with respect to ${\mathcal S}=  k -
\{ {\rm zero}\ {\rm divisors} \}$
we reduce the statement to a semi-local ring $A$.
Then it is well-known (see \cite{Ba}) that
$Pic(A)=0$  and hence the statement follows
from the Hopf-Frobenius case.
\end{proof}

\section{The quantum double of an FH-algebra}

Let $k$ be a commutative ring.
We note that the {\it quantum double} $D(H)$, due to Drinfel'd
\cite{Drin}, is definable for a  finite projective 
Hopf algebra $H$ over $k$:  
 at the level of coalgebras it is given by  
\[
D(H) := H^{*\, {\rm cop}} \otimes_k H ,
\]
where $H^{*\, {\rm cop}}$ is the co-opposite of $H^*$, the coproduct
being $\Delta^{\rm op}$. 

 The multiplication
on $D(H)$ is described in two equivalent ways
 as follows  \cite[Lemma 10.3.11]{Mont}. In terms of  the notation
 $gx$ replacing $g \otimes x$ for every $g \in H^*, x \in H$, both $H$ and $H^{*\, {\rm cop}}$ are subalgebras of
$D(H)$, and for each $g \in H^*$ and $x \in H$,
\begin{equation}
xg := \sum(x_1 g S^{-1}x_3)x_2 
= \sum g_2(S^{-1}g_1 \rightharpoonup x \leftharpoonup g_3).
\end{equation}

The algebra $D(H)$  is a Hopf algebra 
 with antipode $S'(gx) :=  Sx S^{-1}g$, the proof proceeding as in 
\cite{CK2}. A Hopf algebra  
 $H'$ is {\it almost cocommutative},
 if there exists  $R \in H' \otimes H'$, called the
{\it universal $R$-matrix}, such that 
 \begin{equation}
R\Delta(a) R^{-1} = \Delta^{\rm op}(a)
\label{eq:urm}
\end{equation}
 for every $a \in H'$.
A {\it quasi-triangular
Hopf algebra $H'$}  
 is almost cocommutative with universal $R$-matrix  satisfying
  the two equations, 
\begin{eqnarray}
(\Delta \otimes \Id)R & =& R_{13}R_{23}
\label{eq:un} \\
(\Id \otimes \Delta)R &= & R_{13}R_{12}.
\label{eq:deux}
\end{eqnarray}
By a proof like that in \cite[Theorem IX.4.4]{CK2},
$D(H)$ is a quasi-triangular Hopf algebra 
 with universal $R$-matrix
\begin{equation}
R = \sum_i e_i \otimes e^i \in D(H) \otimes D(H),
\label{eq:drin}
\end{equation} 
where $(e_i,e^i)$ is a finite projective base of $H$ \cite{Drin}. 

The next theorem is now a straightforward generalization of
 \cite[Theorem 4.4]{Rad93}.

\begin{theorem}
If $H$ is  an FH-algebra, then the quantum double
$D(H)$ is a unimodular FH-algebra.
\label{th-bfs2}
\end{theorem}
\begin{proof}
  Let $f$ be an FH-homomorphism
with $t$ a right norm. Then  $T := S^{-1}f$ is a  left norm in $H^*$. Let 
$b^{-1}$ be  the left distinguished group-like
 element in $H$ satisfying $Tg = g(b^{-1})T$
for every $g \in H^*$.
Moreover, note that  $\ell := S^{-1}(t)  $ is  a left norm in $H$. 

In this proof we denote elements of $D(H)$ as tensors in $H^* \otimes H$.
We claim that $T \otimes t$ is a left and right integral in $D(H)$.
We first show that it is a right integral. 

The transpose of Formula~\ref{eq:Radford} in Theorem~\ref{prop-gote} is 
 $\sum t_1 \otimes t_2 = \sum t_2 \otimes b^{-1} S^2t_1$.
Applying $ \Delta \otimes S^{-1}$ to both sides yields
$\sum t_1 \otimes t_2 \otimes S^{-1}t_3 =
\sum t_2 \otimes t_3 \otimes (St_1)b$.
It follows easily that 
\begin{equation}
\sum S^{-1}t_3b^{-1} t_1 \otimes t_2 = 
1 \otimes t.
\label{eq:susannah}
\end{equation}

We next make a computation like that in \cite[10.3.12]{Mont}. 
Given a simple tensor $g \otimes x \in D(H)$, 
note that in the second line below we use
 $Tg = g(b^{-1}) T$ for each $g \in H^*$,
and in the third line we use Equation~\ref{eq:susannah}:
\begin{eqnarray*}
(T \otimes t)(g \otimes x) & = & \sum Tg(S^{-1}t_3(-)t_1) 
\otimes t_2x \\
& = &  Tg(S^{-1}t_3b^{-1}t_1) \otimes t_2x \\
& = & g(1) T \otimes tx \\
& = & g(1) \epsilon(x) T \otimes t
\end{eqnarray*}

In order to show that $T \otimes t$ is also a left integral,
we note that 
 Formula~\ref{eq:Radford} applied to the right norm $T' = S^{-1}T$ in
$H^*$ is 
$
\sum T'_1 \otimes T'_2 = 
\sum T'_2 \otimes m^{-1} S^2T'_1
$. Apply $S \otimes S$ to obtain  
\begin{equation}
\sum T_2 \otimes T_1 = \sum T_1 \otimes S^2T_2 m.
\end{equation}
Applying $\Delta  \otimes S^{-1}$ to both sides yields
 $\sum T_2 \otimes T_3 \otimes m S^{-1}T_1 = 
\sum T_1 \otimes T_2 \otimes  S T_3$.
Whence 
\begin{eqnarray}
\sum T_2 \otimes T_3m S^{-1} T_1 & = &
\sum T_1 \otimes T_2 S T_3 \nonumber \\
& = & T \otimes  1.
\label{eq:nanny}
\end{eqnarray}

Then 
\begin{eqnarray*}
(g \otimes x)(T \otimes t) & = & \sum gT_2 \otimes (S^{-1 }T_1 \rightharpoonup
x \leftharpoonup T_3)t \\
& = & \sum gT_2 \otimes
S^{-1}T_1(x_3) T_3(x_1)x_2t \\
& = &  \sum gT_2 \otimes 
[T_3 m S^{-1}T_1](x)t \\
& = & gT \otimes \epsilon(x) t = g(1) \epsilon(x) T \otimes t
\end{eqnarray*}
Thus $T \otimes t$ is also a left integral. 

Next we note that $T \otimes t$ is an FH-homomorphism 
for $D(H)^*$, since $D(H)^* \cong
  H^{\rm op} \otimes H^*$,  the ordinary
tensor product of algebras (recall that 
 $D(H)$ is the  ordinary tensor product of
coalgebras $(H^{\rm op})^* \otimes H$). This follows from 
 $T \otimes t$ being a  right integral in $D(H)$ on the one hand,
while, on the other hand, 
 $H^{\rm op}$ and
$H^*$ are FH-algebras with 
FH-homomorphisms $T = S^{-1}f$ and $t$.

Since $T \otimes t$ is an FH-homomorphism for $D(H)^*$, it follows
that $T \otimes t$ is a right norm in $D(H)$.  Since $T \otimes t$
is a left integral in $D(H)$, it follows that it is a left norm too.
Hence, $D(H)$ is unimodular.
\end{proof}
\begin{remark}
We are thankful to Gigel Militaru for pointing out that the result above
was proved in \cite{CMZ} using a different method.
\end{remark}

\begin{corollary}
  $S(t) \otimes f$ is an FH-homomorphism for $D(H)$.
\end{corollary}
\begin{proof}
Note that $S(t) \otimes f$ is a right integral in 
$D(H)^* \cong H^{\rm op} \otimes H^*$, since $S(t)$ and $f$
are right integrals in $H^{\rm op}$ and $H^*$, respectively. 
Then 
\begin{equation}
(T \otimes t)(S(t) \otimes f) = \epsilon_{D(H)^*} T(S(t))f(t) 
=  \epsilon_{D(H)^*}. 
\label{eq:june10}
\end{equation}
so that $S(t) \otimes f$ is a right 
norm in $D(H)^*$.   By 
Proposition~\ref{cor-june5},  $D(H)$ is an FH-algebra  
with FH-homomorphism
$S(t) \otimes f$. \end{proof}

The next theorem implies that the quantum double
 $D(H)$ of an FH-algebra is a symmetric algebra,
of which   \cite[Corollary 3.12]{BFS2} is a special case.

\begin{theorem}
A unimodular almost commutative FH-algebra $H'$ is a symmetric algebra.
\end{theorem}
\begin{proof}  
Since $H'$ is unimodular, Lemma~\ref{lemma-one}
shows that $H'$ has 
 Nakayama automorphism  $S^2$.  
Since $H'$ is almost commutative, a computation like Drinfeld's  
(cf.\ \cite[Proposition 10.1.4]{Mont})  shows that 
 $S^2$ of an almost commutative Hopf algebra $H$ 
 is an inner automorphism:
if $R = \sum_i z_i \otimes w_i$ is
the universal $R$-matrix satisfying Equation~\ref{eq:urm}, 
then $S^2(a) = uau^{-1}$ where $u = \sum_i (Sw_i)z_i$. Thus, 
the Nakayama automorphism is inner, and $H'$ is a symmetric algebra.
\end{proof}

\vspace{.5cm}

\begin{bf} Acknowledgements. \end{bf}  The  authors thank NorFA of Norway
and NFR of Sweden, respectively, for their support of this paper. 
The first author also thanks Robert Wisbauer and the university in D\"{u}sseldorf 
for their support and hospitality during two recent stays.

\end{document}